\documentclass[12pt]{article}
\usepackage{amsfonts}
\usepackage{epsfig}

\title{A Trichotomy for Rectangles Inscribed in
Jordan Loops}
\author{Richard Evan Schwartz \thanks{\hskip 5 pt Supported by 
N.S.F. Research Grant DMS-1204471}}

\newtheorem{theorem}{Theorem}[section]

\newtheorem{lemma}[theorem]{Lemma}

\newtheorem{corollary}[theorem]{Corollary}

\def\startproof{{\bf {\medskip}{\noindent}Proof: }}

\def\endproof{$\spadesuit$  \newline}

\def\C{\mbox{\boldmath{$C$}}}%
\def\R{\mbox{\boldmath{$R$}}}%
\def\Z{\mbox{\boldmath{$Z$}}}%

\begin{document}
\maketitle

\begin{abstract}
We prove a general structural theorem
about rectangles inscribed in Jordan
loops.  One corollary is that all but
at most $4$ points of any Jordan loop
are vertices of inscribed rectangles.
Another corollary is that a Jordan
loop has an inscribed rectangle of
every aspect ratio provided it has
$3$ points which are not vertices
of inscribed rectangles.
\end{abstract}

\section{Introduction}

A {\it Jordan loop\/} is the image of the circle under a continuous
injective map into the plane.
O. Toeplitz conjectured in 1911 that every Jordan loop contains $4$
points which are the vertices of a square.  This is often called
the {\it Square Peg Conjecture.\/}  An affirmative answer is known
in many special cases.
In 1913,
Emch [{\bf Emch\/}] proved the result for
convex curves. In 1944, L. G. Shnirlmann [{\bf Shn\/}] proved the
result for sufficiently smooth curves.
In 1961, R. Jerrald [{\bf Jer\/}] extended this to the case of
$C^1$ curves. Recently T. Tao [{\bf Ta\/}] proved the result for
special curves having even lower regularity.  The above is a very
partial survey of the literature.
The 2014 survey paper by B. Matschke [{\bf Ma1\/}] and the recent
book by I. Pak [{\bf P\/}] have  extensive discussions
of the history of the Square Peg Conjecture and
many additional references.

There is also some work done in the case of
rectangles. In 1977,
H. Vaughan [{\bf Va\/}] gave a proof that
every Jordan loop has an inscribed rectangle.
A recent paper of C. Hugelmeyer [{\bf H\/}]
combines Vaughan's basic idea  with some
very modern knot theory results to show
that a smooth Jordan loop
always has an inscribed rectangle of
aspect ratio $\sqrt 3$.
The recent paper
[{\bf AA\/}] proves that any quadrilateral inscribed in a
circle
can (up to similarity) be inscribed in any
convex smooth curve.  See also [{\bf Ma2\/}].
In the recent paper [{\bf ACFSST\/}],
the authors show
that every Jordan Loop contains
a dense set of points which are vertices of
inscribed rectangles.
For additional work on
inscribed rectangles, see [{\bf Mak1\/}], [{\bf Mak2\/}],
and [{\bf MW\/}].

Relatedly, one can consider the situation
for triangles. In 1980, M. Meyerson [{\bf M\/}] proved
that all but at most $2$ points of any Jordan loop
are vertices of inscribed equilateral triangles.
This result is sharp because two points of a
suitable isosceles triangle are not vertices
of inscribed equilateral triangles.  In 1992,
M. Neilson [{\bf N\/}] proved that
an arbitrary Jordan loop contains a dense set 
of points which are
vertices of inscribed triangles of any given shape.

We are going to prove a strong version of
Meyerson's Theorem for rectangles.
Let $I(\gamma)$ denote the space of
all labeled rectangles inscribed in $\gamma$.
We always label a rectangle $R$ so that the
vertices of $R$ go counterclockwise around $R$.
We orient $\gamma$ so that it goes counterclockwise
around the region in the plane it bounds.
The space $I(\gamma)$ is naturally a
subset of $\R^8=(\R^2)^4$.

We call a rectangle $R$ inscribed in $\gamma$
{\it graceful\/} if the cyclic order imparted
on the vertices of $R$ by the ordering on
$\gamma$ coincides with the ordering we have
already given to the vertices of $R$.
Let $G(\gamma) \subset I(\gamma)$ denote the space of
gracefully inscribed labeled rectangles.

The {\it aspect ratio\/} of a rectangle in $I(\gamma)$
is the length of the second side divided
by the length of the first side.
Let $\rho(S)$ denote the union of all the
aspect ratios of rectangles in $S$.
Given $S \subset I(\gamma)$, let
$V(\gamma,S) \subset \gamma$ denote the
subset of $\gamma$ consisting of points
which are vertices of rectangles in $S$.

\begin{theorem}[Trichotomy]
\label{threepoint}
Let $\gamma$ be an arbitrary
Jordan loop.  Then $I(\gamma)$
contains a connected set $S$ satisfying
one of the following.
\begin{enumerate}
\item Members of $S$ have uniformly large
area, $V(\gamma,S)=\gamma$, and $1 \in \rho(S)$.
\item Members $S$ have uniformly large
diameter, $\rho(S)=(0,\infty)$, and
$V(\gamma,S)$ contains
all but at most $4$ points of $\gamma$.
\item The set $S$ has members of every
sufficiently small diameter, and
$V(\gamma,S)$ contains all but at most
$2$ points of $\gamma$.
\end{enumerate}
Moreover, $S \subset G(\gamma)$.
\end{theorem}

We mention three corollaries.
The first two corollaries are
immediate, and we will prove
the third one in \S \ref{finalremark}.

\begin{corollary}
  \label{fourpoint}
Let $\gamma$ be any Jordan loop. Then
all but at most $4$ points of $\gamma$
are vertices of rectangles gracefully
inscribed in $\gamma$.
\end{corollary}
This result is sharp: There are $4$ points
of a non-circular ellipse which are
not vertices of any inscribed rectangle.

\begin{corollary}
  \label{twopoint}
Let $\gamma$ be any Jordan loop. If
$\gamma$ has $3$ points which are
not vertices of gracefully inscribed rectangles
then $\gamma$ has gracefully inscribed
rectangles of every aspect ratio.
\end{corollary}

\begin{corollary}
\label{measure}
Let $\gamma$ be any Jordan loop and let
$\mu$ be any non-atomic measure on
$\gamma$ having mass $1$.
Then $\gamma$ has a gracefully inscribed
rectangle $\gamma$ such
that the total $\mu$-measure of each pair
opposite sides of $\gamma$ cut off by $R$
is $1/2$.
\end{corollary}

\noindent
{\bf Remarks:\/}
\newline
(1)  I describe the cases in the
Trichotomy Theorem respectively as {\it elliptic\/},
{\it hyperbolic\/},
and {\it parabolic\/}, because
the geometry of the situation
seems to vaguely resemble the
action of these kinds of
linear  transformations on $\R^2$.
\newline
(2) Note that there are examples, such
as the circle,  for which both the hyperbolic
and elliptic cases occur.
\newline
(3) I conjecture that the parabolic case cannot actually occur.
This conjecture immediately implies the Square Peg Conjecture.
\newline
(4) The elliptic case occurs for
any curve with $4$-fold rotational symmetry
but I conjecture that the hyperbolic case
also occurs for every Jordan loop.
This conjecture implies the much stronger
result that every Jordan loop has an
inscribed rectangle of every aspect ratio -- a
conjecture that is not even known in the smooth or
polygonal cases.
\newline
 (5)
 One could view Corollary \ref{fourpoint}
 as a version
of Meyerson's Theorem [{\bf M\/}] for rectangles, though
our proof is much different.  Our recent preprint
[{\bf S1\/}] gives a proof of Meyerson's Theorem along the
lines of the proof in this paper.
\newline
(6)
The paper [{\bf ACFSST\/}] also has
Corollary \ref{measure} (without the
graceful bit) when
$\gamma$ is rectifiable and
$\mu$ is arc-length measure normalized to have
total length $1$.  The methods are different.
\newline

We prove the Trichotomy Theorem by taking a suitable
limit of the polygonal case. 
Let $\gamma$ be a polygon.
By an {\it arc component\/} of $I(\gamma)$ we
mean a connected component of $I(\gamma)$ which
is homeomorphic to an arc.
By {\it proper\/}, we mean that
as one moves towards
an endpoint of an arc component
in $I(\gamma)$, the aspect ratio
tends either to $0$ or to $\infty$.
Moreover, we insist that the rectangles
at each end of a proper arc accumulate
on a chord of $\gamma$ and that the two
chords are distinct.  The left side of
Figure 1 below suggests an example of a proper arc
in $I(\gamma)$ when $\gamma$ is an equilateral
triangle.  The black segments are the two
chords of accumulation. This arc actually
is a component of $G(\gamma)$.

\begin{theorem}
  \label{generic}
  There is an open dense subset $\cal P$ of
  polygons with the following property.
  For each $\gamma \in \cal P$ the space
  $I(\gamma)$ is
  a piecewise smooth $1$-manifold whose
  arc components are proper.
  Moreover the aspect ratio function
  $\rho: I(\gamma) \to (0,\infty)$ is injective
  in a neighborhood of each smooth point of $I(\gamma)$
  and, $\rho^{-1}(1)$ consists entirely of smooth points.
\end{theorem}

We define $2$ kinds of components of $I(\gamma)$.
\begin{itemize}
\item A component $A$ of $I(\gamma)$
is a {\it  hyperbolic component\/} if
the aspect ratio of the rectangles
in $A$ tends to $0$ as the
rectangles tend towards one
endpoint of $A$, and to $\infty$ as they
tend to the other endpoint of $A$.
The example in Figure 1 is a hyperbolic component.
\item 
The operation of cyclically relabeling
gives a $\Z/4$ action on the space
$I(\gamma)$ which has no fixed points.
We call a  component of $I(\gamma)$
{\it elliptic\/} if it is stabilized by
the $\Z/4$ action.  These components are loops.
\end{itemize}
We call a component
of $I(\gamma)$ {\it global\/} if it
is either hyperbolic or elliptic. The reason for
the name is that $V(\gamma,S)$ contains all but at
most $4$ points of $\gamma$ when $\gamma$ is hyperbolic
(Lemma \ref{extension}) and
all points of $\gamma$ when $\gamma$ is elliptic
(Lemma \ref{coverall}).

The left and right hand sides of Figure 1 respectively
suggest a hyperbolic and an elliptic component.
The figure on the right, drawn by hand, may have
small inaccuracies.  The basic construction
is to have the centers of the rectangles
wind once around a small central curve, while
the rectangles themselves rotate one quarter
of the way around and return to their original
location but with the corresponding cyclic
relabeling.  When this operation is repeated
$4$ times, one has an elliptic component.
It would be simpler to keep the centers fixed,
but I wanted to show an example without $4$ fold
rotational symmetry.

\begin{center}
\resizebox{!}{2.3in}{\includegraphics{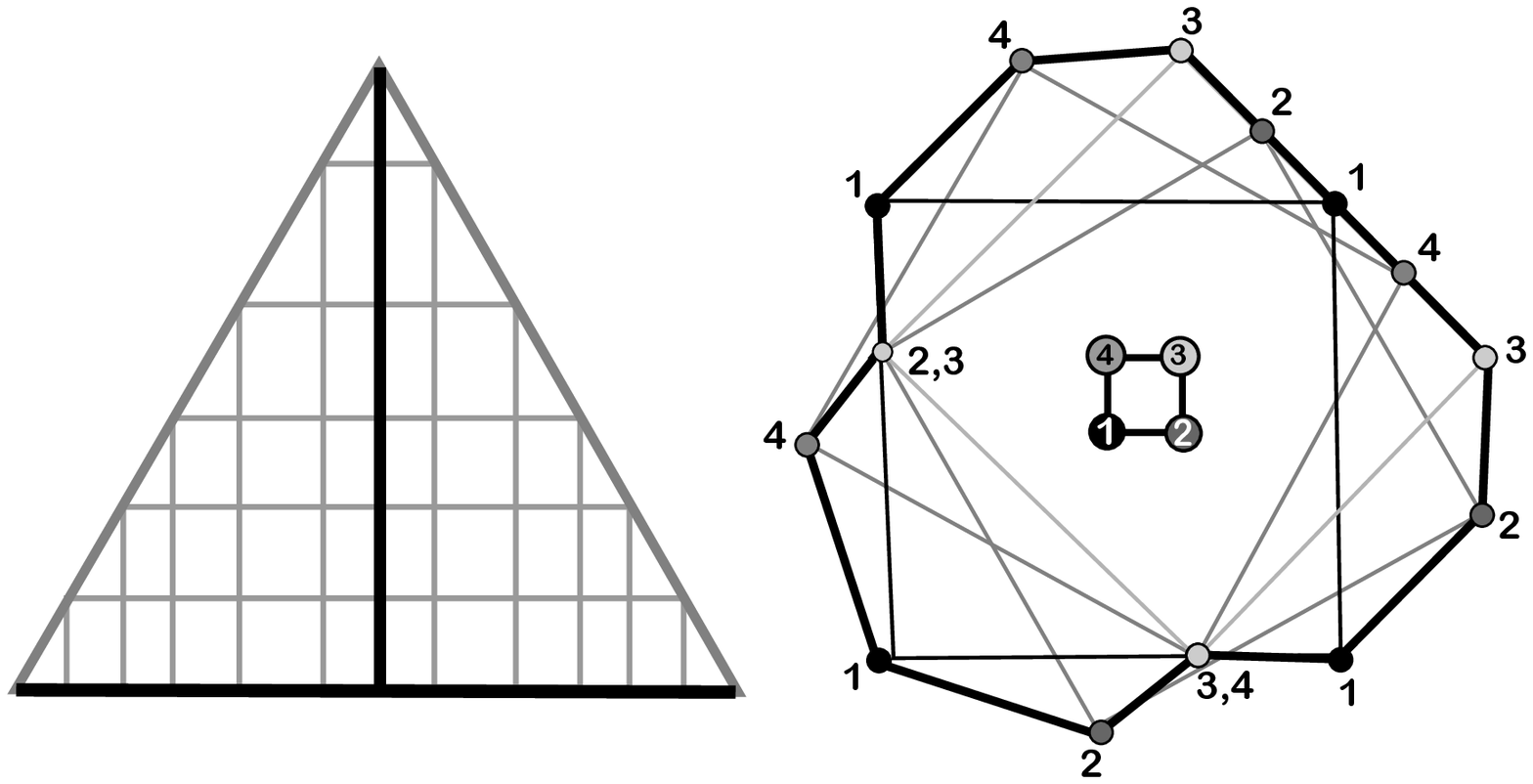}}
\newline
{\bf Figure 1:\/} A hyperbolic component and an elliptic
component.
\end{center}

The relabeling action permutes the various
components of $I(\gamma)$ and we call the
orbits of this action the {\it unlabeled components\/}.
We define the following quantities:
\begin{itemize}
\item  $\Omega(\gamma)$ is the
number of unlabelled inscribed squares.
\item $\Omega_H(\gamma)$ is the number
of unlabelled hyperbolic components.
\item $\Omega_E(\gamma)$
is the number of unlabeled elliptic components.
\end{itemize}
We will establish the following equation for
each $\gamma \in \cal P$.
\begin{equation}
\label{parityX}
\Omega(\gamma)+\Omega_H(\gamma)+\Omega_E(\gamma) \equiv 0\ {\rm mod\/}\ 2.
\end{equation}

It is well known that for the generic polygon the number
of unlabeled inscribed squares is odd.  See for instance
[{\bf St\/}] or [{\bf P\/}, Theorem 23.11].
Hence $I(\gamma)$ always contains a global
component. Finally we prove the following result.
\begin{theorem}
\label{ungood}
For $\gamma \in \cal P$ the only
global components of $I(\gamma)$
belong to $G(\gamma)$.
\end{theorem}
Theorem \ref{ungood} now implies
the following result.

\begin{theorem}
 \label{global}
  For each $\gamma \in \cal P$ the space
  $G(\gamma)$ contains a global component.
\end{theorem}
We get the Trichotomy Theorem by taking
a suitable limit of Theorem \ref{global}.

We also mention another corollary of
Theorem \ref{ungood}:  A generic polygon
has an odd number of gracefully inscribed
squares.  See \S \ref{parity0}.
I don't think that this corollary
follows directly from [{\bf P\/}, Theorem 23.11],
which makes a statement about the parity of
the number of all inscribed squares.
\newline

Here is an outline of the paper.
In \S 2 we will deduce the Trichotomy Theorem
from Theorem \ref{global}.
Following \S 2, the
rest of the paper is about polygons.

In \S 3 we will deduce Equation
\ref{parityX} from Theorem \ref{generic}.

In \S 4 we prove Theorem \ref{generic}.
This is really just an exercise in
transversality, and many methods would work.

In \S 5 we prove Theorem \ref{ungood}.

We warn the reader about one persistent abuse of
terminology. When we
speak of {\it a rectangle in\/} $I(\gamma)$ (or
in related configuration spaces) we mean
the rectangle corresponding to a member of $I(\gamma)$
and not some kind of configuration of $4$ elements
of $I(\gamma)$. We hope that this does not cause confusion.
\newline

One thing I would like to mention is that
I discovered all the results in this
paper by computer experimentation.
I wrote a Java program which
computes the space
$G(\gamma)$ in an efficient way
for polygonal loops $\gamma$ having
up to about $20$ sides.  

I would like to thank Arseniy Akopyan,
Peter Doyle, Cole Hugelmeyer, and Sergei
Tabachnikov for helpful and interesting
conversations related to this paper.
I would also like to thank the referee
of this paper for very helpful comments.
I would like to thank the National Science
Foundation, the Simons Foundation, and
the Isaac Newton Institute for their
generous support while I worked on this paper.

\newpage

\section{The Trichotomy Theorem}

In this chapter we deduce
the Trichotomy Theorem from
Theorem \ref{global} using
a limiting argument.  We begin with
some preliminary material
on point-set topology.

\subsection{Hausdorff Limits}

Suppose that $C$ is a compact metric space.
Let $X_C$ denote the set of closed subsets of
$C$.  We define the {\it Hausdorff distance\/}
between closed $A,B \subset C$ to be
the infimal $\epsilon$ such that
each of the two sets is contained in the
$\epsilon$-tubular neighborhood of the
other one.  This definition makes
$X_C$ into a compact metric space.

\begin{lemma}
  \label{connect}
  Let $\{A_n\}$ be a sequence of nonempty closed
connected subsets of $C$.  Suppose that
this sequence converges to a subset
$A \subset C$ in the Hausdorff metric.
Then $A$ is connected.
\end{lemma}

\startproof
If $A$ is disconnected,
there are disjoint open sets $U,V \subset C$
such that $A \subset U \cup V$, and
$A \cap U$ and $A \cap V$ are both
not empty.  The following pairs of sets are
compact and disjoint:
$$
(A,C-U-V) \hskip 30 pt 
(A\cap U,C-U), \hskip 30pt (A \cap V,C-V).
$$
(The set $A \cap U$ is compact because
$A \cap U=A-V$. Similarly for $A \cap V$.)
Hence, there is some $\epsilon>0$ such that
every point in the first set of a pair is
at least $\epsilon$ from every point in the
second pair.
Therefore, for all sufficiently large $n$, the
set $A_n$ intersects both $U$ and $V$.
Since $A_n$ is connected, this
is only possible if $C-U-V$ contains
a point $x_n \in A_n$.  But then
$x_n$ is at least $\epsilon$ from
$A$, independent of the choice of $n$.
This contradicts the fact
that $A_n \to A$ in the
Hausdorff metric.
\endproof

\noindent
    {\bf Remarks:\/} \newline
    (1) Since $C$ is compact, the set $A$ must be
    non-empty, by the Bolzano-Weierstrass Theorem.
    \newline
    (2) In our application the sets $A_n$ will be
    path connected.  However, 
    there is no guarantee that the limit $A$ is
    path connected as well. Consider a sequence of path approximations
    to the topologist's sine curve.

\subsection{The Circular Invariant}

Let $S^1$ be the unit circle.  Let
$|\alpha|$ denote the arc length of an arc $\alpha \subset S^1$.
Let $\Sigma \subset (S^1)^4$ denote the
subset of distinct labeled quadruples,
which go counterclockwise around $S^1$.
We call these {\it cyclic quadrilaterals\/}.
Let $\sigma_k$ be the $k$th vertex of $\sigma$.
 Any $\sigma \in \Sigma$ defines arcs
$\alpha_0,\alpha_1,\alpha_2,\alpha_3$ with
$\alpha_k \subset S^1-\sigma$ having endpoints
$\sigma_k$ and $\sigma_{k+1}$. The indices are
taken mod $4$. We write $\alpha_k(\sigma)$ when
we want to emphasize the dependence on $\sigma$.

We define the {\it circular invariant\/} 
$\Lambda: \Sigma \to (0,\infty)$ by
\begin{equation}
\label{circular}
\Lambda(\sigma)=\frac{|\alpha_0|+|\alpha_2|}{|\alpha_1|+|\alpha_3|}.
\end{equation}
When $\sigma$ consists of the vertices
of a rectangle, $\Lambda(\sigma)$ is the
aspect ratio of this rectangle.  Otherwise
$\Lambda(\sigma)$ is only vaguely related to
an aspect ratio.

\begin{lemma}
\label{prerect}
Let $\{\sigma_n\}$ be a sequence of cyclic quadrilaterals having
diameter greater than some positive $\delta$ for all $n$, and
circular invariant 
converging to $0$.  Then the arcs
$\alpha_0(\sigma_n)$ and $\alpha_2(\sigma_n)$
shrink to points and
the arcs $\alpha_1(\sigma_n)$ and $\alpha_3(\sigma_n)$ remain
uniformly large.
\end{lemma}

\startproof
The hypotheses imply that 
$\lim_{n \to \infty} |\alpha_0(\sigma_n)|+|\alpha_2(\sigma_n)|=0.$
By the triangle inequality,
$\min(|\alpha_1(\sigma_n)|,|\alpha_3(\sigma_n)|)>\delta/2$
for $n$ large.
\endproof

We call $A \subset \Sigma$ {\it extensive\/} if
$A$ is connected and $\Lambda(A)=(0,\infty)$.

\begin{lemma}
\label{extension}
If $A$ is extensive then all but at most
$4$ points of $S^1$ are vertices of members of $A$.
If, additionally, $A$ contains elements of
arbitrarily small diameter then all but at
most $2$ points of $S^1$ are vertices of
members of $A$.
\end{lemma}

\startproof
Let $\pi_k: A \to S^1$ be the map such that
$\pi_k(\sigma)=\sigma_k$, the $k$th vertex.
The set $J_k=\pi_k(\Sigma)$ is connected, and
therefore either an open arc, a closed arc, a
half-open arc, or all of $S^1$.
If $B=S^1-\bigcup J_k$ contains more than $4$ points,
there is some interval
$B' \subset B$ which has an endpoint in common with
$J_k$ and an endpoint in common with $J_{k+1}$ for some $k$.
But then $B' \subset \alpha_k(\sigma)$ for all $\sigma \in A$.
This bounds $\Lambda(A)$
away from $0$ or $\infty$, depending on the parity of $k$,
a contradiction.

If $B$ has at least $3$ points, then each member of
$A$ is a cyclic quadrilateral which nontrivially
intersects each of $3$ disjoint circular arcs.
There is a uniform positive lower bound to the
diameter of such cyclic quadrilaterals.
\endproof

\subsection{A Compactness Result}

For $K \geq 1$ let
\begin{equation}
\Sigma(K)=\Sigma \cap \Lambda^{-1}([1/K,K]).
\end{equation}
Note that $\Sigma(K)$ is not compact because
for any unit complex number $u$ with positive
imaginary part, the cyclic quadrilateral
$(1,u,-1,\overline u)$ lies in $\Sigma(1)$.
In spite of this problem, we will prove a
compactness result that involves $\Sigma(K)$.

Let $\gamma$ be a general Jordan loop.
We fix, once and for all, some homeomorphism
$\phi: S^1 \to \gamma$.
In [{\bf Tv\/}], Tverberg
gives a way to approximate
$\gamma$ by a sequence
$\{\gamma_n\}$ of parametrized
embedded polygons so that the
parametrizations $\phi_n: S^1 \to \gamma_n$
converge uniformly to $\phi$.
Let $\Gamma$ denote the set of
$\sigma \in \Sigma$ such that
$\phi(\sigma)$ is
the vertex set of a rectangle in $G(\gamma)$,
the space of rectangles gracefully inscribed in $\gamma$.
Likewise define $\Gamma_n$ relative to
$\phi_n$ and $\gamma_n$.

\begin{lemma}
  \label{precompact}
  For each fixed $K \geq 1$ there is a compact
  subset $C(K) \subset \Sigma$ such that
  $\Gamma_n \cap \Sigma(K) \subset C(K)$ for all $n$.
\end{lemma}

\startproof
If this is false then we can pass to a subsequence and
cyclically relabel so that one of the following
two things is true.
\begin{enumerate}
\item There is a sequence $\{\sigma_n\}$, with $\sigma_n \in \Gamma_n \cap \Sigma(K)$,
such that $|\alpha_0(\sigma_n)| \to 0$.
\item There is a sequence $\{\sigma_n'\}$ with $\sigma_n' \in \Gamma_1 \cap \Sigma(K)$,
such that $|\alpha_0(\sigma_n')| \to 0$.
\end{enumerate}

Consider Case 1.  Let
$R_n$ be the rectangle in $G(\gamma_n)$ corresponding to $\sigma_n$
We claim that $|\alpha_1(\sigma_n)| \to 0$ as well.
If not, then
the denominator of the expression for
$A(\sigma_n)$
in Equation \ref{circular} is uniformly
bounded away from $0$.  But then so is
the numerator.  Hence $|\alpha_2(\sigma_n)|$ is uniformly large.
Since $\phi_n \to \phi$ uniformly, the
side of $R_n$ corresponding to $\alpha_0(\sigma_n)$ shrinks to a
point but the opposite side corresponding
to $\alpha_2(\sigma_n)$ does not. This is impossible for
a sequence of rectangles.  This proves our claim.  But now
we can repeat the same argument twice more
to show that $|\alpha_k(\sigma_n)| \to 0$ for $k=0,1,2,3$.
This contradicts the fact $\sum_{k=0}^3 |\alpha_k(\sigma_n)|=2\pi$.

Case 2 is really just an instance of Case 1 relative
to the sequence of polygons $\{\gamma_n'\}$, the
sequence of maps $\{\phi_n'\}$, and the
limit $\phi': S^1 \to \gamma'$. Here (somewhat trivially)
$$\gamma'=\gamma_1'=\gamma_2'=\gamma_3'=... = \gamma_1, \hskip 30 pt
\phi'=\phi_1'=\phi_2'=\phi_3'=... = \phi_1.$$
So, the argument in Case 1 takes care of Case 2.
\endproof

\subsection{Limits of Hyperbolic Components}

We keep the notation from the previous sections.
In this section we will
consider the special case that $G(\gamma_n)$
has a hyperbolic component $H_n$ for all $n$.
Recall that $A \subset \Sigma$ is
extensive if $A$ is
connected and $\Lambda(A)=(0,\infty)$.

\begin{lemma}
  The subset $A_n \subset \Gamma_n$ corresponding
  to $H_n$ is extensive.
\end{lemma}

\startproof
Each end of $H_n$ consists of rectangles which
accumulate on some chord of $\gamma_n$.
Hence, the circular invariants of the
corresponding cyclic quadrilaterals tend
to $0$ or $\infty$, with one case happening
at one end and the other case happening at the other end.
Hence $\Lambda(A_n)=(0,\infty)$.
Also, $A_n$ is homeomorphic to the
arc $H_n$ and hence also connected.
\endproof

\begin{lemma}
  \label{jordan}
  $\Gamma$ contains an extensive subset $A$.
\end{lemma}

\startproof
We keep the same notation as in the previous lemma.
To make our proof more flexible, we only use
the property that $A_n$ is path connected, and contains
a member with circular invariant $n$ and a member with circular invariant $1/n$.
Since $A_n$ is connected, $A_n$ contains some $\sigma_n$ with circular invariant $1$.
For each $K=1,...,n$ we define
$A_n(K)$ to be the minimal arc of $A_n$ which contains
$\sigma_n$ and has endpoints with circular invariant
$1/K$ and $K$ respectively.  By construction $A_n(K)$ is defined for
$n \geq K$ and furthermore $A_n(K) \subset C(K)$,
the compact set from Lemma \ref{precompact}.
Finally,
\begin{equation}
\label{contain0}
 A_n(1) \subset ... \subset A_n(n)
\end{equation}

For fixed $K$, the sequence
$\{A_n(K)\}$ is a sequence of closed connected subsets of
the compact set $C(K)$.  Using Cantor's diagonal trick,
and compactness, we can find a subsequence so that
for each $K=1,2,3,...$ the sequence
$\{A_n(K)\}$ converges to some $A(K) \subset C(K)$ as
$n \to \infty$.
By Lemma \ref{connect} the set $A(K)$ is connected.
Moreover, $A(K)$ contains elements of circular invariant
$K$ and $1/K$.  Each $\sigma \in A(K)$ is such that
$\phi(\sigma)$ is the vertex set of a non-degenerate limit of
gracefully inscribed rectangles.  Hence $\sigma \in \Gamma$.  In short,
$A(K) \subset \Gamma$.
Equation \ref{contain0} gives us $A(2) \subset A(3) \subset A(4)...$.
The nested union of connected sets is connected.
Therefore $A=\bigcup_K A(K)$ is connected.
By construction $A$ is extensive and $A \subset \Gamma$.
\endproof

\subsection{Limits of Elliptic Components}

We continue with the notation above.  This
time we treat the special case where
$G(\gamma)$ has an elliptic component
$E_n$ for all $n$. 
We call the sequence $\{E_n\}$
{\it steady\/} if there is a uniform
positive lower bound to the side length
of any rectangle in any $G(\gamma_n)$,
independent of $n$, and otherwise
{\it wobbly\/}.

\begin{lemma}
  \label{jordan2}
  If $\{E_n\}$ is wobbly
  then $\Gamma$ contains an extensive set $A$.
\end{lemma}

\startproof
Let $A_n \subset \Sigma$ be the set
which corresponds to $E_n$.
Again, $A_n$ is connected.
It cannot be the case that the
circular invariants of members of
$A_n$ are uniformly bounded away
from $0$ and $\infty$.  Otherwise,
the lack of diameter bound
contradicts Lemma \ref{precompact}.
Therefore, after taking a subsequence,
we can arrange that $A_n$ either has
a member with circular invariant $n$ or
a member with circular invariant $1/n$.
Given the invariance of $A_n$
under cyclic relabeling, we see that
$A_n$ has a member with circular
invariant $1/n$ and a member with
circular invariant $n$.
This is all we used in the proof of Lemma \ref{jordan}.
So, the same proof works here as well.
\endproof

\begin{lemma}
  \label{coverall}
  We have $1 \in \rho(E_n)$ and
  $V(E_n,\gamma_n)=\gamma_n$.
\end{lemma}

\startproof
Recall that $\rho$ is the aspect ratio.
Since $E_n$ is invariant under cyclic relabeling,
we have that $r \in \rho(E_n)$ iff $1/r \in \rho(E_n)$.
Since $E_n$ is connected, $1 \in \rho(E_n)$.
This is the first claim.

Let $v_k(R)$ denote the $k$th vertex of a rectangle $R$.
We take indices mod $4$.
Choose any rectangle $R_0 \in E_n$.
Since $E_n$ is an elliptic component, there is a path
$\{R_t|\ t \in [0,1]\}$ of rectangles in $E_n$
such that $v_k(R_t)$  connects
$v_k(R_0)$ and $v_{k+1}(R_0)$.
We write $v_k(t)=v_k(R_t)$.
It suffices to prove $\gamma_n=\bigcup v_k([0,1])$.

For the proof, we identify
$\gamma_n$ with $\R/4\Z$ so that the vertices of
$R_0$ are $[0],[1],[2],[3]$. 
The path $v_k$ connects $[k]$ to $[k+1]$.
Let $\widehat v_k: [0,1] \to \R^2$ be the lift of
$v_k$ such that $\widehat v_k(0)=k$.
Note that $\widehat v_k(0) \not = \widehat v_k(1) \in \Z$.
Hence the interval $I_k=[\widehat v_k(0),\widehat v_k(1)]$
has length at least $1$. 
We have
\begin{equation}
\label{ineq00}
\widehat v_0(t)<\widehat v_1(t)<\widehat v_2(t)<\widehat v_3(t)<\widehat v_1(t)+4.
\hskip 15 pt \forall t \in [0,1].
\end{equation}
Equation \ref{ineq00} holds for $t=0$ and
any failure at time $t$ 
would result in the points $\{v_k(t)\}$
not being distinct. Hence $I_{k+1}=I_k+1$.
This immediately implies that $\bigcup I_k$
contains an interval of length $4$ and
hence so does $\bigcup \widehat v_k([0,1])$.
Hence $\R/4\Z \subset \bigcup v_k([0,1])$.
\endproof

\begin{lemma}
  \label{steady}
  If $\{E_n\}$ is steady then
  $G(\gamma)$ contains a compact connected
  subset $S$ such that $1 \in \rho(S)$ and
$V(\gamma,S)=\gamma$.
\end{lemma}

\startproof
Let $A_n \subset \Sigma$ be the
subset corresponding to $E_n$.
By hypotheses there is a single
compact subset $C \subset \Sigma$
such that $A_n \subset C$ for all $n$.
Passing to a subsequence, we take
the Hausdorff limit $A=\lim A_n$.
The set $A$ is connected, by Lemma \ref{connect}.
We let $S=\phi(A)$. For the same reason as in
the hyperbolic case, $S \subset G(\gamma)$.

Since $1 \in \rho(E_n)$ for all $n$,
we can find some square in $S$ corresponding
to a limit of uniformly large squares in
$E_n$. Hence $1 \in \rho(S)$. 

Let $p \in \gamma$ be any point.
Let $p_n \in \gamma_n$ be such that
$p_n \to p$.  Since $V(E_n,\gamma_n)=\gamma_n$
we can find a rectangle $R_n \in E_n$ such
that $p_n$ is a vertex of $R_n$.
There is a uniform lower bound to the
side lengths of these rectangles.
Hence, any limit $\lim R_n$ will be
a rectangle in $S$ having $p$ as
a vertex. Hence $V(S,\gamma)=\gamma$.
\endproof

\subsection{The Main Argument}

Perturbing our polygons if necessary,
we can assume that each $\gamma_n$ satisfies
Theorem \ref{global}.
Passing to a subsequence, we reduce to
either the hyperbolic case considered
above, the steady elliptic case, or
the wobbly elliptic case.
\newline
\newline
{\bf Case 1:\/}
In the steady elliptic case,
Lemma \ref{steady}
gives Option 1 of the Trichotomy Theorem.
\newline
\newline
{\bf Case 2:\/} In
the hyperbolic case or the wobbly elliptic
case, Lemma \ref{jordan} or Lemma \ref{jordan2}
guarantees that $\Gamma$ contains an extensive set $A$.
Let $S \subset G(\gamma)$ be the subset corresponding to $A$.
Suppose that there
is a uniform positive lower bound to the diameters of members of $A$.
We show that Option 2 of the Trichotomy Theorem holds.

By Lemma \ref{extension}, at most $4$ points of
$S^1$ are not vertices of members of $A$.
Hence $V(S,\gamma)$ contains all
but at most $4$ points of $\gamma$.

Since $S$ is connected,
$\rho(S)$ is connected.
Since $A$ is extensive, $A$ contains a sequence
$\{\sigma_n\}$ of circular quadrilaterals
whose circular invariant tends to $0$.
By Lemma \ref{prerect}, this is only possible if
the two arcs $\alpha_0(\sigma_n)$ and
$\alpha_2(\sigma_n)$
in Equation \ref{circular} shrink to points
and the other two arcs remain uniformly long.
But then the aspect ratios of the corresponding
rectangles tend to $0$.  Hence
$\rho(S)$ contains points arbitrarily near $0$.
The same argument shows that $\rho(S)$ contains
points arbitrarily near $\infty$.
Hence $\rho(S)=(0,\infty)$.
\newline
\newline
{\bf Case 3:\/}
The only remaining case is that $\Gamma$ contains an
extensive set $A$ without the positive lower diameter bound.
In this case, Lemma \ref{extension} shows that
at most $2$ points of $S^1$ are not vertices of
members of $A$.  Hence $V(S,\gamma)$ contains all
but at most $2$ points of $\gamma$.  Since $A$
contains members of every sufficiently small
diameter, the set $S$ does as well.
This gives us Option 3 of the Trichotomy Lemma.

\subsection{Non-Atomic Measures}
\label{finalremark}

Here we prove Corollary \ref{measure}.
Suppose that $\mu$ is a non-atomic probability
measure on $\gamma$.
Call a quadrilateral gracefully inscribed in
$\gamma$ {\it nice\/} if it cuts $\gamma$
in such a way that opposite arcs have
$\mu$-measure $1/2$.  We are looking
for a nice inscribed rectangle.
Since $\mu$ is non-atomic, we can choose
our homeomorphism $\phi$ so that $\phi$
pushes forward arc length on $S^1$ to
$2\pi \mu$.  Then nice rectangles correspond
to elements of $\Gamma$ having circular 
invariant $1$. 

In Cases 2 and 3 of our proof above, the
sets $A$ are extensive and have such
cyclic quadrilaterals. So, Corollary
\ref{measure} is true in Cases 2 and 3 above.

For Case 1, we revisit the proof of Lemma \ref{steady}.
Since $A_n$ is invariant under cyclic
relabeling, we have 
$1 \in \Lambda(A_n)$.  So, by 
Lemma \ref{precompact}
we can take a limit and get $1 \in \Lambda(A)$.
The corrsponding cyclic quadrilateral in $A$
corresponds to a nice rectangle in the set $S$.

\newpage

\section{The Parity Equation}

\subsection{Outline of Proof}

In this chapter we deduce Equation \ref{parityX}
from Theorem \ref{generic}.
Let $\cal P$ be the open dense set of polygons
from Theorem \ref{generic}.  We fix some
$\gamma \in \cal P$ for the entire argument.
The space $I(\gamma)$ of labeled inscribed rectangles
is a $1$-manifold, by Theorem \ref{generic}.
The cyclic group $\Z/4$ acts on $I(\gamma)$ by
cyclically relabeling the rectangles.  Again,
the labeling of a rectangle goes counterclockwise
around the the rectangle.  This is a free action:
No point of $I(\gamma)$ is fixed by the relabeling.

For emphasis,
we call the components of
$I(\gamma)$ {\it labeled\/}. 
An {\it unlabeled\/} component is the orbit of a
labeled component under the labeling action.
We define
the {\it order\/} of a labeled component
 to be the number
of labeled components in its orbit -- either
$1$, $2$, or $4$.

We say that a labeled
rectangle $R$ is {\it associated\/} with
the labeled component that contains the point
representing $R$.  We say that a labeled
rectangle $R$ is associated to an orbit
of a labeled component if it is associated
to one of the labeled components in the
orbit.  Finally, we say that an unlabeled
rectangle is associated to an unlabeled
component if the corresponding labeled
rectangles are associated with the 
corresponding orbit.

Below we will prove the following $4$ claims.
\begin{enumerate}
\item The number of unlabeled inscribed squares
associated to an unlabeled hyperbolic component
is odd.  
\item The number of unlabeled inscribed squares
associated to any other unlabeled arc component is even.
\item The number of unlabeled inscribed squares
associated to an unlabeled elliptic component
is odd.
\item The number of unlabeled inscribed
squares associated to any other unlabeled
loop component is even.
\end{enumerate}
Combining these claims, we see that
the total number of unlabeled inscribed
squares, namely $\Omega$, has the same
parity as $\Omega_H+\Omega_E$, the
total number of unlabeled hyperbolic
components plus the total number of
unlabeled elliptic components.
This is exactly
Equation \ref{parityX}.

\subsection{The Arc Components}

\begin{lemma}
Every arc component of $I(\gamma)$ has order $4$.
\end{lemma}

\startproof
Recall that $I(\gamma)$ is naturally a
subset of $\R^8$.  Given an arc component
$\zeta$, there are two points $\zeta_1,\zeta_2$
in $\R^8$ corresponding to the ends of $\zeta$.
Up to cyclic relabeling, each of these points has
the form
$(a,b,a,b,c,d,c,d)$
where $(a,b) \not = (c,d)$.
Each of these points encodes the 
chord of $\gamma$ corresponding to an
end of $\zeta$, and $\zeta_1 \not = \zeta_2$.

Suppose that $\psi(\zeta)=\zeta$ for
some cyclic relabeling map $\psi$.
Given the form of our points, we see that
that $\zeta_j \not = \psi(\zeta_j)$.
This means that $\psi$ must interchange
$\zeta_1$ and $\zeta_2$.
But this means that $\psi$ is a homeomorphism of
the arc $\zeta$ which swaps its ends.  This
situation forces $\psi$ to fix a point of
$\zeta$.  This is impossible.
\endproof

\begin{lemma}
Each labeled hyperbolic component contains an odd number of
labeled inscribed squares and any other labeled
arc component contains an even number of labeled
inscribed squares.
\end{lemma}

\startproof
Let $\zeta$ be a labeled hyperbolic arc component.
Let $\rho: \zeta \to (0,\infty)$ be the aspect
ratio function.  By definition
$\rho(p)=1$ if and only if $p$ represents a
square. At one end of $\zeta$, the value of
$\rho$ is less than $1$.  At the other end, the
value of $\rho$ is greater than $1$.  Given
that $\rho$ is injective in a neighborhood of
each point of $\rho^{-1}(1)$, this means
that $\rho=1$ an odd number of times on $\zeta$.
The argument for the non-hyperbolic arcs is the same
except that
$\rho$ is either greater than $1$ at both
ends of $\zeta$ or less than $1$ at both ends.
\endproof

\noindent
{\bf Proof of Claim 1:\/}
Let $\zeta_1,\zeta_2,\zeta_3,\zeta_4$ be the
hyperbolic components comprising an orbit.
There is some odd $k$ such that $\zeta_1$ has
$k$ labeled inscribed squares associated to it.
But then, by symmetry, the same holds for
the other components.  Hence, there are a total
of $4k$ labeled inscribed squares associated to
these components.  But this means that there
are $k$ unlabeled inscribed squares associated to
these components. \endproof
\newline
\newline
{\bf Proof of Claim 2:\/}
The proof of Claim 2 is the same as the proof of
Claim 1 except that now $k$ is even.
\endproof

\subsection{The Loop Components}

\begin{lemma}
A labeled elliptic component contains $4k$
labeled inscribed squares for some odd integer $k$.
\end{lemma}

\startproof
Let $\zeta$ be some labeled elliptic component and
$r_0$ be the point of $\zeta$ that has aspect ratio
less than $1$.
Let $r_1,r_2,r_3$ be the successive images of $r_0$ under the
relabeling map.   Let
$\zeta_k$ be the arc of $\zeta$ bounded by $r_k$ and $r_{k+1}$
with indices taken mod $4$.  Consider the restriction
of the aspect ratio function $\rho$ to $\zeta_0$.
We have $\rho(r_1)=1/\rho(r_0)$.  So, as we trace
out $\zeta_0$ from $r_0$ to $r_1$ we see that
$\rho$ starts out less than $1$ and ends up
greater than $1$.  Hence $\rho$ attains the value
$1$ an odd number $k$ of times on $\zeta_1$.
By symmetry, $\rho$ attains the value $1$ exactly
$k$ times on each arc $\zeta_k$. This give
a total of $4k$.
\endproof

\noindent
{\bf Proof of Claim 3:\/}
Let $\zeta$ be some labeled elliptic component.
The orbit of $\zeta$ is just $\zeta$ itself.
We have just seen that the number of
labeled inscribed squares associated to
$\zeta$ is $4k$ for some odd $k$.
But then the number of unlabeled
inscribed squares associated to $\zeta$ is
$k$. \endproof

\begin{lemma}
A labeled loop component contains $2k$
labeled inscribed squares for some integer $k$.
If the component has order $2$ then $k$ is even.
\end{lemma}

\startproof
Let $\zeta$ be such a component.
Since $\zeta$ is a topological loop and
$\rho$ is injective in a neighborhood
of each point of $\rho^{-1}(1)$, the
map $\rho$ attains the value $1$ an
even number of times.  This is the
first statement.

Suppose then that $\zeta$ has order $2$.
Let $r_0$ be some point of $\zeta$ such
that $\rho(r_0) \not = 1$.  Since
$\zeta$ has order $1$, the element of
$\Z/4$ which sents vertex $0$ to vertex $2$
must be the one which stabilizes $\zeta$.
Let $r_1$ be the image of $r_0$ under
this relabeling element.  Both $r_0$
and $r_1$ have the same aspect ratio.
Hence $\rho$ attains the value $1$ an even
number of times on each of the arcs
of $\zeta$ joining $r_0$ to $r_1$.
\endproof

\noindent
{\bf Proof of Claim 4:\/}
Let $\zeta$ be a labeled loop
component which is not elliptic.
Regardless of whether $\zeta$
has order $1$ or $2$, the
preceding lemma says that there
are $8h$ labeled squares 
associated to the orbit of $\zeta$.
Hence there are an even number of
unlabeled squares associated to the
orbit of $\zeta$. \endproof

\subsection{Discussion}
\label{parity0}

The analysis above combines with
Theorem \ref{ungood} to prove the
follownig result.

\begin{corollary}
\label{oddgrace}
For $\gamma$ in $\cal P$ the space
$G(\gamma)$ contains an odd number of squares.
\end{corollary}

\startproof
The analysis above shows that
each non-global component
of $I(\gamma)$ 
contributes an even number to the total
number of unlabeled inscribed squares.
By Theorem \ref{ungood}, therefore,
$I(\gamma)-G(\gamma)$ contains an
even number of inscribed squares.
Hence $G(\gamma)$ contains an
odd number of inscribed squares.
\endproof

Here is a way to deduce Theorem \ref{global}
from Corollary \ref{oddgrace} without
appealing to Theorem \ref{ungood}.
The same argument as above establishes
Equation \ref{parityX} when we make
the counts with respect to components
in $G(\gamma)$.  So, if we already know
that there are an odd number of
gracefully inscribed squares, then
our new count implies $G(\gamma)$ has
either a hyperbolic component or an
elliptic component.  This approach to
Theorem \ref{global} is more direct,
and we originally took it.  However,
our direct proof that $G(\gamma)$ has
an odd number of squares, a homotopy
argument, was rather tedious.

\newpage

\section{The Moduli Space}

\subsection{Inscribing Rectangles in Four Lines}
\label{main}

For this chapter we work in
the complex plane $\C$.  
Let $L=(L_0,L_1,L_2,L_3)$ be a quadruple of lines in
$\C$.  We assume throughout the chapter that these
lines are in general position.
We say that a rectangle $R$ 
{\it is gracefully inscribed in\/} $L$ if the
vertices $(R_0,R_1,R_2,R_3)$ go
cyclically around $R$ (either clockwise
or counterclockwise) and satisfy
$R_i \in L_i$ for $i=0,1,2,3$.
We let $G(L)$ denote the set of
rectangles gracefully inscribed in $L$.
We think of $G(L)$ as a subset of
$\R^8$.  
In [{\bf S2\/}] I worked out quite a bit
of the structure of $G(L)$.  Here I will give a
more abstract and less detailed treatment,
but when relevant I will point out
the stronger results that appear in
[{\bf S2\/}].

We define the {\it aspect ratio\/}
$\rho: G(L) \to \R$ by the formula
\begin{equation}
\rho(R)=\pm \frac{|R_2-R_1|}{|R_1-R_0|}.
\end{equation}
The sign is $-1$ if $R$ is
clockwise ordered and $+1$ if
$R$ is counterclockwise ordered.
We allow $\rho$ to be both positive
and negative, though ultimately
we just care about the case $\rho>0$.

\begin{lemma}
\label{intersect}
Let $G(L,\rho)$ denote the subset of
$G(L)$ consisting of rectangles having
aspect ratio $\rho$. For generic choice of $L$,
the space $G(L,\rho)$ has at most one element
for every choice of $\rho$.
\end{lemma}

\startproof
Let $R$ be a rectangle in $G(L)$. We have
$$R_2-R_1=i\rho(R_1-R_0), \hskip 30 pt
R_3-R_2=R_1-R_0.$$
Writing this as a matrix equaion, we have
$(R_2,R_3)=M(R_0,R_1)$,
where
\begin{equation}
M=\left[\matrix{-i\rho & 1+i\rho \cr 1-i\rho & i\rho}\right].
\end{equation}
Note that $\det(M)=1$ and ${\rm trace\/}(M)=0$.
This means that $M$ is always invertible and
indeed an involution.    Let $\Pi_{12}=L_0 \times L_1$
and $\Pi_{34}=L_2 \times L_3$. These are both totally
real planes in $\C^2$.  The solutions we seek are
the points of $M(\Pi_{12}) \cap \Pi_{34}$.
These two planes are either disjoint, or
intersect in a $d$-dimensional affine subspace
for $d=0,1,2$. 

The case $d=2$ is certainly not generic, and we rule
out the case $d=1$ with a dimension count.
The space of rectangles of aspect ratio $\rho$ is $5$ dimensional.
So, the space of pairs of rectangles having aspect $\rho$ is
$9=5+5-1$ dimensional.   Given two such rectangles
$R$ and $R'$, we can recover the quadruple $L$ by letting
$L_k$ be the line through $R_k$ and $R_k'$.  This accounts
for all quadruples of interest.  The space of
rectangles of aspect ratio $\rho$ gracefully inscribed in $L$ is $1$-dimensional
affine subspace.  Thus, we have overcounted the quadruples
of lines of interest by $2$ dimensions:  Any pair of
rectangles in the family would produce the same quadruple.
This the space of quadruples containing infintely many
rectangles of the same aspect ratio has dimension $7$.
On the other hand, the space of quadruples of lines
has dimension $8$.
\endproof

\noindent
{\bf Remark:\/}
In [{\bf S2\/}] I show that $G(L)$ contains infinitely
many rectangles of the same aspect ratio if and only
if the line through $L_0\cap L_1$ and $L_2 \cap L_3$
is perpendicular to the line through $L_1 \cap L_2$ and $L_3 \cap L_0$.
\newline

Since $L$ is generic, $G(L)$ contains exactly one
rectangle $R_{\rho}$ of aspect ratio $\rho$ provided
that it contains any.
Let $\rho(L)$ denote the set of aspect ratios of
rectangles in $G(L)$.
We define the point $\phi_k(\rho) \in L_k$ denote 
the $k$th vertex of $R_{\rho}$. This gives us a
map $\phi_k: \rho(L) \to L_k$.  We call
$\phi_k$ a {\it vertex map\/}.

Let $\cal L$ be the space of generic quadruples.
We can identity this space with an
open subset of $\R^8$.  We can think of
our vertex map $\phi_k$ as a map
${\cal L\/} \times \R \to \C$.
The point $\phi_k(L,\rho)$ is the vertex of
the rectangle $R_{\rho}$ defined relative
to the configuration $L$.  The domain for
$\rho_k$ is naturally the fiber bundle
$\cal E$ over $\cal L$ whose fiber is $\rho(L)$.

\begin{lemma}
The space $\cal E$ is an open subset of $\R^9$ and
each $\rho_k$ is an analytic function on $\cal E$.
\end{lemma}

\startproof
Referring to Lemma \ref{intersect}, the existence
of $R_{\rho}$ means that the two planes
$M(\Pi_{12})$ and $\Pi_{34})$ are transverse.
A small change in $\rho$ or in $L$ does not change that fact.
Hence $\cal E$ is open in $\R^9$.  The desired point
can be found using linear algebra with inputs that
vary analytically with the coordinates on $\cal E$.
Everything in sight is algebraic, and hence
analytic.
\endproof

\noindent
{\bf Remark:\/} In [{\bf S2\/}] I show that
the set of centers of rectangles in $G(L)$
is a hyperbola minus $2$ points, namely
those corresponding to degenerate rectangles
of aspect ratio $0$ and $\infty$.
There are two unequal values $a_1$ and $a_2$,
corresponding to the points at infinity of
the hyperbola, such that $\rho(L)=\R-\{0,a_1,a_2\}$.
Moreover $a_1a_2$ equals the cross ratio of
the slopes of the lines of $L$.
\newline

There is one degenerate case we need to consider.
We say that a {\it repeating quadruple\/} is
one of the form
\begin{equation}
(L_0,L_0,L_1,L_2), \hskip 15 pt
(L_0,L_1,L_1,L_2), \hskip 15 pt
(L_0,L_1,L_2,L_2), \hskip 15 pt
(L_0,L_1,L_2,L_0).
\end{equation}
Here $L_0,L_1,L_2$ are distinct and
non-parallel lines.
We call these kinds of quadruples
{\it repeating quadruples\/}.
We make the same definition for $G(L)$.
The same results as above apply in this
case.  Indeed, it never happens that
$G(L)$ contains infinitely many rectangles
of the same aspect ratio.
\newline
\newline
{\bf Permutation Trick:\/}
We have now defined two kinds of quadruples,
the generic ones and the generic repeating
ones.  So far we have been talking about
gracefully inscribed rectangles.  Since we are interested
in all inscribed rectangles, rather than just
the gracefully inscribed ones, we note that the permutation
of a generic quadruple is still generic and
the cyclic or dihedral permutation of a
repeating quadruple is still a repeating
quadruple.  For instance, if we are 
interested in the space $G'(L)$ of rectangles
having the property that 
$R_0 \in L_1$ and $R_1 \in L_3$ and
$R_1 \in L_0$ and $R_3 \in L_2$
then we are really considering the space
$G(L')$, where $L'=(L_1,L_3,L_0,L_2)$ is
the suitable permutation of the lines of $L$.
The space $I(\gamma)$ divides into different
subspaces, depending on the combinatorics of
the labelings.  The space $G(\gamma)$ is one
of these.  When proving things for $I(\gamma)$,
we will often specialize to the case of
$G(\gamma)$, with the understanding that the
permutation trick just discussed promotes
the proof we give for $G(\gamma)$ to a
proof for $I(\gamma)$.
\newline
\newline
{\bf Regular and Singular Values:\/}
Since $\phi_k$ is analytic, there are
finitely many values $b_1,...,b_{\ell} \in \rho(L)$
such that $d\phi_k/da=0$.  Here we are differentiating
with respect to the aspect ratio parameter.
We call $b_1,...,b_{\ell}$ the
{\it singular ratios\/} for $\phi_k$ and we call
their images on $L_k$ the {\it singular images\/}.
We call $a \in \R$ a {\it regular ratio\/} if
it is not a singular ratio.  Geometrically,
the $k$th vertex of $R_a$ varies monotonically
on $L_k$ near a regular value $a$. We call
the image of a regular ratio a {\it regular image\/}.

\subsection{The Polygon Space}

In this section, we define our space $\cal P$ of polygons 
and then show that it is open and dense in the
space of all polygons.  What really mean is
that there is a space ${\cal P\/}(N)$ for
each $N \geq 3$, which is a subset of the
finite dimensional space of all labeled $N$-gons,
and ${\cal P\/}(N)$ is open and dense in this space.
The space of all labeled $N$-gons is simply
$\R^{2N}$.  So, our space ${\cal P\/}(N)$ is
an open dense subset of $\R^{2N}$.

Let $\gamma$ be a polygon, with sides
$E_1,...,E_N$.  We insist first of all that
no two sides of $\gamma$ are parallel.
This assumption implies that any rectangle
gracefully inscribed in $\gamma$ has its vertices in at least
$3$ distinct edges.
 We say that an
{\it associated quad\/} $L=(L_0,L_1,L_2,L_3)$ is
a quadruple of lines, cyclically ordered (counterclockwise)
and extending some sides of $\gamma$, which is
either repeating or ordinary in the sense above.
We say that $\gamma$ belongs to $\cal P$ if
\begin{enumerate}
\item All associated quads are generic or repeating.
\item No inscribed rectangle has more than $1$ vertex in common with $\gamma$.
\item No inscribed square has a vertex in common with $\gamma$.
\item No vertex of $\gamma$ is a singular image with
respect to an associated quad.

\end{enumerate}

\begin{lemma}
$\cal P$ is open and dense in the space of all polygons.
\end{lemma}

\startproof
Condition 1 clearly holds on an open
dense set.  The remaining conditions
are open because their negation is
closed. For instance, if we have a
convergent sequence of $N$-gons
having an inscribed square that
shares a vertex with the polygons,
then we can take a limit and get
such a square on the limiting polygon.
The other conditions are similar.

We will deal with Condition 2.
Given any pair $v_1,v_2$ of vertices of
$\gamma$, let $e$ be the edge joining
$v_1$ and $v_2$ and let $\lambda_1$ and
$\lambda_2$ be the lines perpendicular
to $e$ through $v_1$ and $v_2$ respectively.
Let $\lambda_3$ be the circle having
$e$ as a diameter.  We can perturb so that
$\lambda_i \cap \gamma$ is always a finite
set of points.  We can further perturb so that
there are no points
$v_i \in \lambda_i \cap \gamma$ for
$i=1,2$ or $v_1,v_2 \in \gamma_3$ such that
$\|v_1-e_1\|=\|v_2-e_2\|$.
Thus, after finitely many steps, we
get Condition 2 and keep Condition 1.

Conditions 3 and 4 involve a single
vertex at a time.  To show density
for these, we introduce a {\it slide move\/}.
This is defined relative to a vertex $v$ and
one of the edges $e$ incident to $v$.
We replace $v$ by a vertex $v' \in e$
very close to $v$ and then consider the
new polygon having $v'$ as a vertex in
place of $v$ and all other vertices
the same.

If we have a polygon which fails to
have one of the conditions above, we
can associate that failure to a
triple $(v,e,L)$ where $v$ is an
involved vertex, $e$ is an involved
edge, and $L$ is an associated quadruple
one of whose sides extends $e$.  In
case the same problem -- e.g. a square
sharing a vertex with the polygon -- involves
more than one triple $(v,e',L')$ we count
this as a separate problem.  Each slide
move, if done with respect to a sufficiently
nearby vertex, removes one of the problems
and does not create any new ones. 

For instance, given $(v,e,L)$, there
exist points $v' \in e$ arbitrarily
close to $v$ such such that the
square in $G(L)$ does not contain $v'$
and $v'$ is a regular value for the
relevant vertex map.  This follows
from the analyticity of everything
in sight: The problem points on
$e$ are isolated. 

So, we go around making small perturbations
fixing one problem at a time until we are
done, and we can make these perturbations
as small as we like.
\endproof

\subsection{The Manifold Structure}

Now we prove Theorem \ref{generic}.
Let $\gamma$ be a polygon in $\cal P$.

\begin{lemma}
\label{manifold2}
The space $I(\gamma)$ is a piecewise smooth manifold.
\end{lemma}

\startproof
We will prove this for $G(\gamma)$.
As discussed above, the permutation trick
promotes the proof to a proof for
$I(\gamma)$.
We have a partition
\begin{equation}
\label{manifold}
G(\gamma)=G_0(\gamma) \cup G_1(\gamma).
\end{equation}
The points of $G_k(\gamma)$ correspond to
gracefully inscribed rectangles which have $k$ points in common
with the vertex set of $\gamma$.

Each point of $G_0(\gamma)$ corresponds to a
rectangle of $G(L)$ for some unique associated
quadruple of lines.  All nearby points of $G_0(\gamma)$ are
associated to the same quadruple of lines.  Thus $G_0(\gamma)$
is open in $G(\gamma)$, and every point of
$G_0(\gamma)$ has a neighborhood which is
just a copy of a neighborhood of the corresponding
point of $G(L)$.  So, by the results in
\S \ref{main}, the set $G_0(L)$ is a smooth manifold
and the aspect ratio function $\rho$ gives a
coordinate chart.

Let $p \in G_1(L)$.  Let $R$ be the associated
rectangle and let $v \in R$ be the vertex of $R$
which is also a vertex of $\gamma$.  There are
exactly $2$ associated quadruples
$L$ and $L'$ such that the rectangle $R$ associated
to $p$ lies in $G(L)$ and $G(L')$. 
After cyclically relabelling, we can arrange
that $L_0$ and $L_0'$ are the two lines
extending the edges of $\gamma$ incident to
$v$, and $L_j'=L_j$ for $j=1,2,3$.

Let $U$ and $U'$ denote small open subsets of
$p$ in $G(L)$ and $G(L')$ respectively.
Each member of $G(\gamma)$ sufficiently close
to $p$ lies in one of $G(L)$ or $G(L')$,
so a small neighborhood of $p$ in $G(\gamma)$
is given by 
$$(U \cap G(\gamma)) \cup (U' \cap G(\gamma)).$$
Because $p$ is a regular image with respect
to $L$ or $L'$, the set $U$ intersects
$G(\gamma)$ in a half-open interval having
$p$ as endpoint.  The idea here is that as
we vary the point in $G(L)$ the vertex near
$p$ move monotonically along $L_0$,
spending half the time on the edge of
$\gamma$ contained in $L_0$ and half the
on $L_0$ outside this edge. The same
goes for $U'$.  Thus,
the two half-open arcs fit together
to give a neighborhood of $p$ in $G(\gamma)$
homeomorphic to an arc.

We have shown that every point of
$G(\gamma)$ has an arc neighborhood,
and every point of $G_0(\gamma)$ is
smooth. Finally, it follows from
the analyticity of the vertex maps that
there are only finitely many rectangles
of $G(\gamma)$ having any given point of
$\gamma$ as a vertex. In particular, $G_1(\gamma)$ is
just a finite set of points. Hence
$G(\gamma)$ is a piecewise smooth $1$-manifold.
\endproof

\begin{lemma}
The aspect ratio
function $\rho$ is locally injective at
each smooth point of $I(\gamma)$, and
$\rho^{-1}(1)$ consists entirely of
smooth points.
\end{lemma}

\startproof
The set of smooth points is
precisely the set $G_1(\gamma)$ considered
above.  The restriction of
$\rho$ to a small neighborhood of such a
point coincides with the restriction
of $\rho$ to some neighborhood of the
corresponding point of $I(L)$.
This restriction is injective by
Lemma \ref{intersect}.
The second
statement of the lemma is exactly
Condition 3.
\endproof

\begin{lemma}
Each arc component of $I(\gamma)$ is proper.
\end{lemma}

\startproof
We prove this for $G(\gamma)$ and then
use the permutation trick to promote
the proof to one for all of $G(\gamma)$.
We first recall what we are trying to prove.
We have $G(\gamma) \subset \R^8$.
Let $A \subset G(\gamma)$ be an arc component.
Let $\partial A=\overline A-A$ denote the boundary
of $A$ in $\R^8$.  We will show that
$\partial A$ consists of $2$ distinct points,
both representing chords of $\gamma$.

In the proof of Lemma \ref{manifold2} we
saw that the space $G_1(\gamma)$ is a
finite set of points.  So, if a rectangle
in $G(\gamma)$ has sufficiently large or
small aspect ratio it must be a smooth
point.  This means that there are two
unique associated quadruples $L$ and $L'$
such that the ends of $A$ respectively
lie in $G(L)$ and $G(L')$. 
If the rectangles at one end of $A$
accumulate to something other than
a line segment, then we can take
a subsequential limit of these
uniformly large and fat rectangles
to get another member of $G(L)$.
This rectangle would also grace $L$,
and hence $\gamma$,
and it would have a neighborhood in
$G(\gamma)$ that
overlaps with $A$. This is a contradiction.
Hence, as we exit an end of $A$ the
corresponding rectangles accumulate on
a line segment.

For the end of $A$ associated to $L$,
the relevant rectangles accumulate
to a chord that has one vertex on
two consecutive lines of $L$ and one
vertex on the other two.  This
chord is uniquely determined by $L$
and, by general position, uniquely
determines $L$ amonst all associated quads.
All the same remarks apply to the
other end of $A$, which is associated
to the quadruple $L'$.

It remains to show that our chords
are distinct.  If not, then $L=L'$
and the two ends, which both have their
elements in $G(L)=G(L')$, have
some rectangles in commn. This is a
contradiction.
\endproof

\noindent
{\bf Remark:\/} 
The ends of each proper arc of $G(\gamma)$ are
critical points for the distance function
$d: \gamma \times \gamma \to [0,\infty)$,
at least after one makes a suitable definition
for what this means in the polygonal case.
After doing thousands of experiments, I
noticed that one end of a proper arc is
always a saddle point (i.e. a point with
Morse index $1$) and the other end is
always either a maximum or a minimum
(i.e. a point with Morse index $0$ or $2$.)
I have no idea how to prove it, but this
fact suggests hidden depths.

\newpage

\section{Global Components are Graceful}

\subsection{The Easy Part}
\label{easy}

In this chapter we prove Theorem \ref{ungood}.
We say that a rectangle $R \in I(\gamma)$ is
{\it ungracefully inscribed\/} $\gamma$ if the
cyclic order imparted on $R$ from
the ordering on $\gamma$ is the clockwise
ordering of the vertices of $R$.
Let $G^*(\gamma)$ denote the space
of rectangles which are
ungracefully inscribed in $\gamma$.

\begin{lemma}
\label{ungood1}
\label{ungood2}
A global component of $I(\gamma)$ lies
in $G(\gamma)$ or $G^*(\gamma)$,
\end{lemma}

\startproof
Consider the hyperbolic case first.
If the rectangles in the hyperbolic
component $A$ are neither
gracefully nor ungracefully inscribed, then 
there are a pair of opposite sides
of the rectangles such that the
endpoints on one pair of opposite
sides interlace on $\gamma$ with
the endpoints on the other pair.
However, at one end of $A$, these
edges are very short.  This is
only possible if the rectangles are
shrinking to a single point.
This is a contradiction.

When $A$ is elliptic, we revisit
the proof of Lemma \ref{coverall}.
The same paths $\{v_k\}$ exist in
the more general setting, but now
there is some permutation
$\{i_0,i_1,i_2,i_3\}$ of $\{0,1,2,3\}$
such that $v_k(1)=[i_k]$.  Equation
\ref{ineq00} holds in this setting,
and tells us that there are
integers $j_0,j_1,j_2,j_3$ such that
$$i_0+j_0<i_1+j_1<i_2+j_2<i_3+j_3<i_0+j_0+4.$$
Here $i_k+j_k=\widehat v_k(1)$.
This forces $[i_0],[i_1],[i_2],[i_4]$ to be
consecutive residue classes in $\Z/4$.
But then $(i_0,i_1,i_2,i_3)$ is a
cyclic permutation of $(0,1,2,3)$.
This happens if and only if
$A \in G(\gamma)$ or $A \in G^*(\gamma)$.
\endproof

\noindent
{\bf Remark:\/}
The components $G(\gamma)$ and $G^*(\gamma)$ might look
superficially similar, but actually they are quite different.
For instance, $G^*(\gamma)$ is empty if $\gamma$ is convex.
Before we prove that $G^*(\gamma)$ has no global components,
we explain why most of the Trichotomy
Theorem follows from what we have already done.
The definition of the circular invariant and the
other proofs in \S 2, go through practically word for
word if we work with components in $G^*(\gamma)$ rather
than in $G(\gamma)$. 
Thus, if we use Lemma \ref{ungood1} in place of
Theorem \ref{ungood}, we get the whole Trichotomy
Theorem except that the last statment is 
weaker: $S \subset G(\gamma)$ or $S \subset G^*(\gamma)$.

\subsection{The Elliptic Case}

In this section we prove that
$G^*(\gamma)$ contains no elliptic
components.  We will deduce this result
from a theorem about inscribed triangles. 
We define gracefully and ungracefully
inscribed triangles the same way as for rectangles.
We say that an {\it essential graceful loop\/}
(respectively {\it essential ungraceful loop\/}) is
a continuous loop of gracefully  (respectively ungracefully)
inscribed triangles such that each vertex winds a nontrivial
number of times around $\gamma$.  We prove the following result.

\begin{theorem} 
\label{antigrace}
No polygon has an ungraceful loop. \end{theorem}
If we had an elliptic component in $G^*(\gamma)$ we could look at
the loop of triangles made from the first $3$ points.  
The same lifting argument as in
the proof of Lemma \ref{ungood2} shows
that the $k$th vertex of the rectangle family
winds around $\gamma$ a nonzero number
of times.  So, we would get an essential
ungraceful loop,
contradicting Theorem \ref{antigrace}.
We prove Theorem \ref{antigrace} through
two lemmas.

\begin{lemma}
\label{grace1}
A polygon arbitrarily close to $\gamma$
in the Hausdorff metric supports a
graceful essential loop.
\end{lemma}

\startproof
We describe a motion of the points $a,b,c$.
We can find an arbitrarily nearby polygon whose
convex hull has $8$ consecutive vertices
$v_1,...v_8$ which agree with the
vertices of a regular polygon.
We start with points $a,b,c$ located at $v_5,v_6,v_7$.
We then move $c$ all the way around $\gamma$ counterclockwise
until $c=v_4$.  Next, we move $b$ around to $v_3$,
then $a$ around to $v_2$.
Now we have $a,b,c$ located at $v_2,v_3,v_4$.
Finally, we slide this triangle over to 
its original location.
\endproof

If an elliptic component of  $H(\gamma)$
lies in $G^*(\gamma)$ then $\gamma$ supports
an essential ungraceful loop.  But then
so do all polygons sufficiently near
$\gamma$. (The new loops won't necessarily be comprised
of right triangles, but we don't care.)
But then we could have an example of a triangle
which supports both a graceful essential
loop and an ungraceful essential loop.
So, the following lemma finishes the
proof of Theorem \ref{antigrace}.

\begin{lemma}
\label{exclusion}
No polygon can support both an essential graceful
loop and an essential ungraceful loop.
\end{lemma}

\startproof
Let $\Omega$ denote the subset of
distinct triples of points in $\gamma$,
with the order induced by the ordering
on $\gamma$.  Let $\Omega_+ \subset \Omega$
be the subset corresponding to triangles
having positive signed area.  Likewise
define $\Omega_-$.
Our graceful and ungraceful essential
loops respectively define essential
loops $\beta_+$ and $\beta_-$ in
$\Omega_+$ and $\Omega_-$.  We can replace
$\beta_+$ and $\beta_-$ by nearby polygonal loops.

There are nonzero integers $n_{\pm}$ such that
$n_+ \beta_+$ and $n_- \beta_-$ are
homologous.  But then we can find a
piecewise linear surface-with-boundary that has
$n_+ \beta_+$ and $n_- \beta_-$ as
boundary components.  The common
boundary of $\Delta_+$ and $\Delta_-$
is piecewise algebraic, and so
(after we perturb to put things in
general position) the 
intersection $\Delta_+ \cap \Sigma$
consists of finitely many piecewise smooth
loops.  The union of these loops is homologous
to $n_+ \beta_+$, so some component $\beta_0$
is essential.

Say that a {\it stick\/} is a triple of
collinear points of $\gamma$. The loop
$\beta_0$ corrsponds to a continuous family
of sticks having the property that each point
of the stick winds a nonzero number of
times around $\gamma$.  But then there will
be a moment when the middle point of the
stick will be a vertex of the convex hull
of $\gamma$. At this moment, the other
two points of the stick cannot lie on
$\gamma$, and we have a contradiction.
\endproof

\subsection{The Hyperbolic Case}

Now we prove that $G^*(\gamma)$ has
no hyperbolic component.
The argument is similar.
 This time let $\Omega$ denote
the set of quadruples of distinct points in
$\gamma$, with the order induced by the ordering
on $\gamma$.   The space $\Omega$ is homeomorphic
to the product of a $3$-ball and a circle.
 Let $\Omega_+$ (respectively $\Omega_-$) denote the set of quadruples 
defining embedded quadrilaterals, not necessarily convex, whose
ordering goes counterclockwise (respectively clockwise)
around their boundary.  
Unlike in the elliptic case, the two sets $\Omega_+$
and $\Omega_-$ do not partition $\Omega$. That
does not bother us.

We distinguish $2$ special subsets of $\partial \Omega$.
Let $\partial_0 \Omega$ be the set of embedded quadruples
$(a,a,b,b)$ and let $\partial_1 \Omega$ be the
set of quadruples $(a,b,b,a)$.  A hyperbolic
component defines a curve joining
$\partial_0 \Omega$ to $\partial_1 \Omega$.  The
interior of this curve lies in $\Omega_+$ when
the hyperbolic component is graceful and in
$\Omega_-$ when the hyperbolic component is
ungraceful.

\begin{lemma}
\label{grace3}
A polygon arbitrarily close to $\gamma$
in the Hausdorff metric supports a path
in $\Omega_+$ joining $\partial_0 \Omega$
to $\partial_1 \Omega$.
\end{lemma}

\startproof
We describe a motion of points $a,b,c,d$.
Make the same modification to $\gamma$ as in
the proof of Lemma \ref{grace1}.  
(We really don't need all $8$ points here.)
Start out with $a=b=v_3$ and $c=d=v_6$.
First move $b$ and $c$ respectively to
$v_4$ and $v_5$.  Now move $d$ all the
way around $\gamma$ to $v_2$.  At every
stage of this construction we have a
convex quadrilateral.  Now we have
$a,b,c,d$ at $v_2,v_3,v_4,v_5$. Finally,
move $b$ and $c$ back to $v_2$ and $v_3$
respectively.  Now we have
$b=c=v_2$ and $a=d=v_5$.
\endproof

\begin{lemma}
There cannot be paths both in
$\Omega_+$ and $\Omega_-$ joining
$\partial_0 \Omega$ to $\partial_1 \Omega$.
\end{lemma}

\startproof
The proof here is similar to the elliptic
case.  We first perturb so that both paths
are polygonal.  We then observe that
both paths represent the generator of
the relative homology group 
$$H_1(\Omega,(\partial_0(\Omega) \cup \partial_1(\Omega)) =\Z.$$
(The space in question is homotopy equivalent relative the
boundary to an annulus relative its boundary.)
So, we can build a piecewise linear surface-with-boundary
$\Sigma$ whose boundary is made up of a path
$\partial_0 \Omega$, a path in 
$\partial_1 \Omega$, and our $2$ curves.
After we put things in general position, the
intersection $\Sigma \cap \partial \Omega_+$ 
is a union of loops and arcs in $\Sigma$,
one of which joins $\partial_0 \Omega$ to
$\partial_1 \Omega$.  Call this path $\beta_0$.

The quadrilaterals in $\beta$ have $4$ distinct
points but are not embedded.  One possibility is
that the points are all collinear and the other
possibility is shown in Figure 3.

\begin{center}
\resizebox{!}{.5in}{\includegraphics{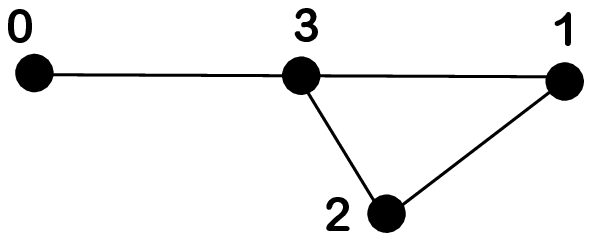}}
\newline
{\bf Figure 3:\/} A degenerate quadrilateral
\end{center}

In case the points are not all collinear, there
is always one segment
that has another vertex between it. In Figure 3,
the segment is $v_0v_1$ and the middle
point is either $v_2$ or $v_3$.  We call
$v_0,v_1$ a {\it framing segment\/}.
We can always cyclically relabel so that at some
point along $\beta_0$ the segment $v_0v_1$ is a
framing segment.  Perturbing $\beta$ slightly,
we can arrange that the points along $\beta$
corresponding to $4$ collinear points are isolated.
Such totally collinear
configurations have codimension $1$ in $\partial \Omega_+$.

Consider what happens to a configuration
with $v_0v_1$ as a framing segment
as we pass through a totally collinear
configuration. Figure 4 shows the only
$3$ possibilities.

\begin{center}
\resizebox{!}{2.6in}{\includegraphics{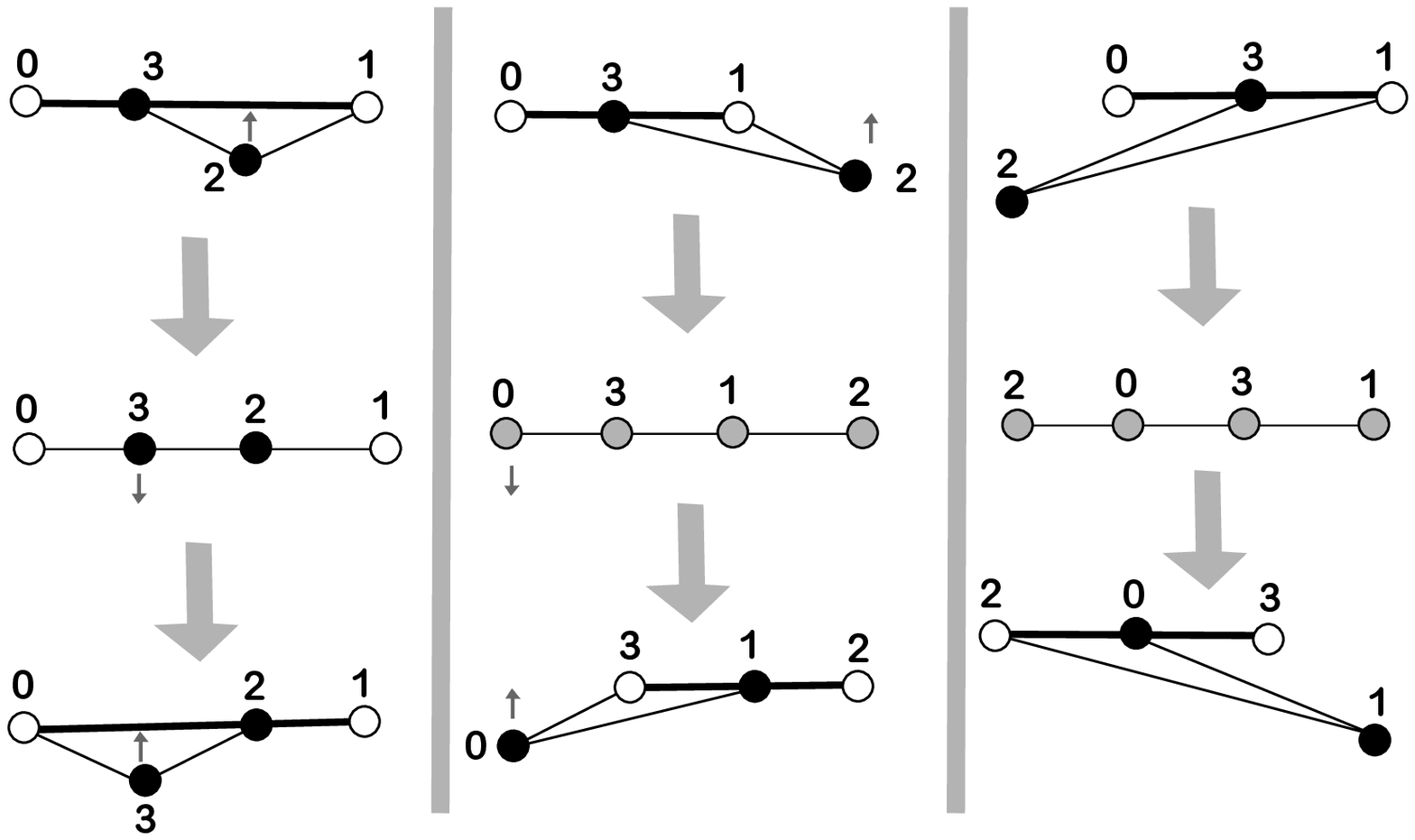}}
\newline
{\bf Figure 4:\/} The only allowable transitions
\end{center}

The framing segment either remains
$v_0v_1$ or else changes to $v_2v_3$.
Therefore, one of the two segments
$v_0v_1$ or $v_2v_3$ is the framing
segment at each point of $\beta_0$.
At $\beta_0 \cap \partial_0 \Omega$,
we have $v_0=v_1$ and $v_2=v_3$.
But the distance between the point
on the framing segment to either
endpoint of the draming segment
therefore tends to $0$.  This shows
that there is a sequence of points
in $\Omega$ converging to a point of
$\partial_0 \Omega$ such that
the diameter of some $3$ of the points
tends to $0$.  This contradicts the
fact that the points of $\partial_0 \Omega$
come together in pairs and not in triples.
\endproof

\newpage

\section{References}
\noindent
[{\bf AA\/}] A. Akopyan and S Avvakumov, {\it Any cyclic quadrilateral can be
inscribed in any closed convex smooth curve.\/} 
arXiv: 1712.10205v1 (2017)
\newline
\newline
[{\bf ACFSST\/}] J. Aslam, S. Chen, F. Frick, S. Saloff-Coste,
L. Setiabrate, H. Thomas, {\it Splitting Loops and necklaces:
Variants of the Square Peg Problem\/},  arXiv 1806.02484 (2018)
\newline
\newline
 [{\bf CDM\/}], J. Cantarella, E. Denne, and J. McCleary,
{\it transversality in Configuration Spaces and the
Square Peg Problem\/}, arXiv 1402.6174 (2014).
\newline
\newline
[{\bf Emch\/}] A. Emch, {\it Some properties of closed convex curves in the plane\/},
Amer. J. Math. {\bf 35\/} (1913) pp 407-412.
\newline
\newline
[{\bf H\/}] C. Hugelmeyer, 
{\it Every Smooth Jordan Curve has an inscribed
rectangle with aspect ratio equal to $\sqrt 3$.\/}
arXiv 1803:07417 (2018)
\newline
\newline
[{\bf Jer\/}]. R. Jerrard, {\it Inscribed squares in plane curves\/},
T.A.M.S. {\bf 98\/} pp 234-241 (1961)
\newline
\newline
[{\bf Mak1\/}] V. Makeev, {\it On quadrangles inscribed in a closed curve\/},
Math. Notes {\bf 57(1-2)\/} (1995) pp. 91-93
\newline
\newline
[{\bf Mak2\/}] V. Makeev, {\it On quadrangles inscribed in a closed curve and
vertices of the curve\/}, J. Math. Sci. {\bf 131\/}(1) (2005) pp 5395-5400
\newline
\newline
[{\bf Ma1\/}] B. Matschke, {\it A survey on the
Square Peg Problem\/},  Notices of the A.M.S.
{\bf Vol 61.4\/}, April 2014, pp 346-351.
\newline
\newline
[{\bf Ma2\/}] B. Matschke, {\it  Quadrilaterals inscribed in
convex curves\/}, \newline
arXiv 1801:01945v2
\newline
 \newline
[{\bf M\/}] M. Meyerson {\it Equilateral Triangles and
Continuous Curves\/}, Fundamenta Mathematicae, 110.1, 1980, pp 1-9.
 \newline
 \newline
[{\bf N\/}] M. Neilson, {\it Triangles Inscribed in Simple Closed
Curves\/}, Geometriae Dedicata (1991)
 \newline
 \newline
[{\bf NW\/}] M. Neilson and S. E. Wright, {\it Rectangles inscribed in
symmetric continua\/}, Geometriae Dedicata {\bf 56(3)\/} (1995) pp. 285-297
\newline
\newline
[{\bf S1\/}] R. E. Schwartz, {\it On Spaces of Inscribed Triangles\/},
preprint 2018.
\newline
\newline
[{\bf S2\/}] R. E. Schwartz, {\it Four lines and a rectangle\/},
preprint 2018.
\newline
\newline
[{\bf Shn\/}], L. G. Shnirelman, {\it On certain geometric properties of closed curves\/}
\newline
(in Russian), Uspehi Matem. Nauk {\bf 10\/} (1944) pp 34-44; 
\newline
available at
http://tinyurl.com/28gsys.
\newline
\newline
 [{\bf St\/}], W. Stromquist, {\it Inscribed squares and square-like
quadrilaterals in closed curves\/}, Mathematika {\bf 36\/} (1989) pp 187-197
\newline
\newline
[{\bf Ta\/}], T. Tao, {\it An integration approach
to the Toeplitz square peg conjecture\/}
\newline
Foum of Mathematics, Sigma, 5 (2017)
\newline
\newline
[{\bf Tv\/}], H. Tverberg, 
{\it A Proof of the Jordan Curve Theorem\/},
\newline
Bulletin of the London Math Society, 1980, pp 34-38.
\newline
\newline
[{\bf Va\/}], H. Vaughan, {\it Rectangles and simple closed curves\/},
Lecture, Univ. of Illinois at Urbana-Champagne.

\end{document}